\documentclass{amsart}
\usepackage{mathrsfs} 
\usepackage{mathrsfs}
\usepackage{tikz}
\usepackage{hyperref}
\newcommand\LP{\textsf{LatticePolytopes}\ }
\newcommand\Mac{\textsf{Macaulay2}\ }
\renewcommand\L{\mathscr{L}}

\newcommand\R{\mathbb{R}}
\newcommand\Z{\mathbb{Z}}
\renewcommand\P{\mathbb{P}}
\renewcommand\O{\mathscr{O}}

\title[Lattice Polytopes in Macaulay2]{LatticePolytopes: A package for computations with lattice polytopes in Macaulay2}
\author{Anders \ Lundman} 
\address{Department of Mathematics, KTH Royal Institute of Techonlogy, SE-100 44 Stockholm, Sweden}
\thanks{The first author is supported by the Swedish Reseach Council, grant number 2014-4763.}
\email{alundman@kth.se}
\author{Gustav \ S\ae d\'{e}n St\aa hl}
\thanks{The second author is supported by the Swedish Reaserch Council, grant number 2011-5599.}
\address{Department of Mathematics, KTH Royal Institute of Techonlogy, SE-100 44 Stockholm, Sweden}
\email{gss@math.kth.se}
\begin{document}
\maketitle

\begin{abstract}
We introduce the package LatticePolytopes for Macaulay2. The package provides methods for computations related to Cayley structures, local positivity and smoothness for lattice polytopes. 

\end{abstract}

\section{Introduction}
A defining characteristic of toric geometry is the interplay between algebraic and convex geometry. For example, affine toric varieties correspond to polyhedral cones, and complete embeddings of projective toric varieties correspond to lattice polytopes. This is the theoretical basis for the two packages  \textsf{NormalToricVarieties} and \textsf{Polyhedra} \cite{Polyh} in \textsf{Macaulay2} \cite{M2}. These packages make it possible to study the core objects and properties of this classical theory. In the present note we introduce the package \textsf{LatticePolytopes}, which extends the functionality of the two above mentioned packages. Our package provides methods for investigating properties of toric varieties relating to Cayley structures, Gauss maps, local positivity, adjunction theory and smoothness, all of which are active research areas, see for example \cite{primer},\cite{PolyAdj},\cite{Lazar}. 

As the name suggests, the package \LP  deals primarily with lattice polytopes. If $M\cong \Z^n$ is a lattice, then a convex lattice polytope $P\subset M\otimes_\R \R=\R^n$ is called smooth if the edge-directions at every vertex form a basis of $M$. Such polytopes are important since they, by the toric dictionary, correspond to smooth polarized toric varieties $(X,\L)$, see \cite[Theorem~2.4.3]{Cox}. Furthermore, the set $P\cap M$ of lattice points contained in $P$ corresponds to a basis of the global sections of $\L$. Using this correspondence, Bogart et al.\ proved that, for a fixed positive integer $N$, there are, up to isomorphism, only finitely many complete embeddings of smooth toric varieties in $\P^N$ \cite{finiteSLP}. Equivalently, in convex geometric terms, there are only finitely many smooth lattice polytopes with $N$ lattice points up to isomorphism. Motivated by this result, all smooth two dimensional lattice polytopes with at most $12$ lattice points, as well as all three dimensional lattice polytopes with at most $16$ lattice points, have been 
classified in \cite{lorenz} and \cite{myclass} respectively. Both these classifications are implemented  in the package \textsf{LatticePolytopes} via the functions \textsf{listSmooth2D} and \textsf{listSmooth3D}.

As can been seen from these classifications, most complete embeddings ${X\hookrightarrow \P^N}$ of smooth toric surfaces and threefolds with small enough $N$ are projective bundles over a smooth toric base. The associated smooth polytopes are in turn endowed with a (strict) Cayley structure \cite{SandraLinFib}. 
Our package provides a method for investigating if a given, not necessarily smooth, rational polytope has a Cayley structure via the boolean method \textsf{isCayley}. Moreover, the method \textsf{Cayley} provides an easy method  to construct such polytopes. 

Another part of the package \LP deals with local positivity. Loosily speaking, local positivity can be described as the study of positive line bundles on a variety. The positivity of a line bundle is given by various generalizations of what it means for a line bundle to be ample. 
Two such possible generalisations are jet seperations and Seshadri constants. 
In the package \textsf{LatticePolytopes}, the methods \textsf{isJetSpanned} and \textsf{degreeOfJetSeperation} are included in order to help the user investigate the jet seperation of line bundles on toric varieties. Furthermore, the package includes the method \textsf{epsilonBounds}, which is a numerical algorithm to bound Seshadri contants at general points based on bounds in \cite[Theorem~3.6]{degen}.

The last part of the package deals with Gauss maps on toric varieties. A Gauss map is a rational morphism which outside the singular locus assigns to every point the projective tangent spaces at that point. Gauss maps have been studied quite intensely in both algebraic and differential geometry. They are also related to local positivity via a generalization to linear spaces of higher order tangency \cite{C22, Pohl,GH79}. For toric varieties, there exist combinatorial descriptions of these maps, see \cite{FukukawaIto} for the classical case and \cite{HGM} in the case of higher order Gauss maps. Unfortunately, the combinatorics quickly becomes complicated, as the number of lattice points of the associated lattice polytopes grows. This can be remedied using our package  both in the classical and higher order case. For the classical case one can use the methods \textsf{gaussFiber} and \textsf{gaussImage} while the methods  \textsf{gausskFiber}‚  and \textsf{gausskImage} apply to the higher order case.

\section{Examples}
Polyhedral sets are implemented in \Mac as objects of the class \textsf{Polyhedron}.
A \textsf{Polyhedron} is a type of hashtables defined within the package \textsf{Polyhedra}. The same package allows for manipulation of general polyhedral sets whereas the methods introduced in our package are unique for lattice polytopes. Using the package \textsf{LatticePolytopes}, the following set of commands creates a list of all smooth lattice polytopes of dimension $2$ with at most 12 lattice points, and then counts how many of them have exactly $n$ lattice points, for each $n\le 12$. 
\begin{verbatim}
L=listSmooth2D();
tally(apply(L,x->#latticePoints(x)))
 = Tally{3 => 1, 4 => 1, 5 => 1, 6 => 3, 7 => 3, 8 => 4, 
        9 => 5, 10 => 7, 11 => 6, 12 => 10}
\end{verbatim}
For instance, we see that there are precisely $3$ smooth lattice 2-polytopes with exactly $7$ lattice points. By the toric dictionary, these numbers correspond to the number of embeddings of smooth toric surfaces $X\hookrightarrow \P^n$ up to isomorphism for every $n\le 11$. 

An interesting characteristic of the set of smooth 3-dimensional lattice polytopes with at most 16 lattice points, classified in \cite{myclass}, is that the vast majority has a so called Cayley structure. The following commands generate this list of polytopes, and then check how many of them are of this type using the method \textsf{isCayley}.
\begin{verbatim}
L=listSmooth3D();
tally(apply(L,isCayley))
= Tally{false => 4, true => 99}
\end{verbatim}
We see that only $4$ of the $103$ polytopes lack a Cayley structure. In general, if the polytope $P$ associated to a polarized toric variety $(X,\L)$ has a Cayley structure, then there is a canonical and explicitly given birational morphism $X'\dashrightarrow X$, where $X'$ is a projective bundle over a smooth toric basis \cite{SandraLinFib}. Moreover, if $X$ is smooth and of dimension at most $3$, then $P$ has a Cayley structure if and only if $X$ is a projective bundle over a smooth toric base \cite[Proposition ~2]{localpos}.

One natural operation on a smooth variety is taking blow-ups in closed subvarieties. Blowing up an embedded toric variety in a torus invariant subvariety is equivalent to cutting off faces of the associated polytope. Here, cutting off a face $Q$ of a polytope $P$ means intersecting $P$ with a certain halfspace not containing $Q$. By the above mentioned equivalence, this is the same as taking the blow-up of the variety corresponding to $P$ at the subvariety corresponding to $Q$. 
These toric blow-ups can be computed in \LP using the command \textsf{toricBlowUp}. The command also asks the user to specify a positive integer $k$ and returns the polytope corresponding to the blow-up $\pi\colon\widetilde{X} \to X$, with exceptional divisor $E$, embedded in projective space via $\pi^*\L-kE$.  As an example, the following lines of code compute the polytope corresponding to the blow-up of $\P^1\times \P^1$ at a fixpoint (corresponding to a vertex of the polytope) embedded by $\pi^*\O(2,2)-E$.
\begin{verbatim}
P=hypercube 2;
BP=toricBlowUp(P,convexHull(matrix{(vertices P)_0}),1)
{vertices P, vertices BP}
 =  {| -1 1  -1 1 |, | -1 0  1  -1 1 |}
     | -1 -1 1  1 |  | 0  -1 -1 1  1 |
\end{verbatim}
%

The package \LP also includes a family of methods to investigate the local positivity properties of polarized toric varieties in terms of the associated polytopes. For example, the command \textsf{degreeOfJetSeparation} computes the largest integer $k$ such that a given line bundle $\L$ on a toric variety is $k$-jet spanned at a specified point. Recall that if $x\in X$ is a smooth point of a projective variety, with maximal ideal $\mathfrak{m}_x$, then a line bundle $\L$ on $X$ is said to be $k$-jet spanned at $x$ if the natural map 
\[j^k\colon H^0(X,\L)\rightarrow H^0(X,\L\otimes \O_{X}/\mathfrak{m}_x^{k+1})\]
 is onto. In local coordinates around a point $x$, this map is given by sending a section to the first $k$ terms of its Taylor expansion at $x$. Thus, jet separation naturally generalize the very ampleness, where the associated embedding separates points and tangent directions. 

Another measure of local positivity is the Seshadri constant at a point $x$ of a smooth polarized variety $(X,\L)$. It is defined as 
\[\epsilon(X,\L;x)=\inf_{C\subseteq X} \frac{C\cdot \L}{m_C(x)},\] 
where the infimum is taken over all irreducible curves $C\subset X$ passing through $x$ and $m_C(x)$ denotes the multiplicity of $C$ at $x$. These constants are closely related to the Seshadri condition for ampleness and the Nagata conjecture. They are in general very hard to compute. For smooth toric varieties, the Seshadri constant at a torus fixpoint coincides with the degree of jet separation at that point \cite{primer}. Upper and lower bounds of Seshadri constants at general points on toric surfaces can be computed in the package \LP using the command \textsf{epsilonBounds}. These combinatorial bounds are based on the bounds introduced in \cite{degen}, and a description of the algorithm appears in \cite[Algorithm~2]{localpos}. Applying these methods to $(\P^1\times \P^1,\O(2,2))$ works as follows.
\begin{verbatim}
P=hypercube 2;
degreeOfJetSeparation(latticePoints(P),transpose(matrix{{1,1}}))
=  2
epsilonBounds(P,17)
 = {2, 2}
\end{verbatim}
Here, the command \textsf{degreeOfJetSeparation} computes the degree of jet separation of $(\P^1 \times \P^1,\O(2,2))$ at the point which is the image of the identity $(1,1)$ of the torus action. Recall that the  degree of jet separation is constant on the big torus orbit of a toric variety and that the general points of the variety lie in this orbit \cite{Perkinson}. The last command shows that $\epsilon(\P^1\times \P^1,\O(2,2);1)$ is bounded from above and below by $2$. Thus, by the above computation, we conclude that the degree of jet separation of the line bundle $\O(2,2)$ on $\P^1\times \P^1$ at the general point, and the Seshadri constant at such a point, coincide.
%

The last part of the package \LP provides methods to understand the fiber and image of the Gauss map and its higher order analogues on toric varieties. For an $n$-dimensional projective variety $X\hookrightarrow \P^N$ which is $k$-jet spanned at the general point, the Gauss map of order $k$ is defined as the rational map $X\dashrightarrow Gr\bigl(\binom{n+k}{k}-1,N\bigr)$ that sends a point $x\in X$ to the $k$-th osculating space of the embedded variety at $x$. Recall that the projectivisation of the image of the map $j^k$, defining $k$-jet spannedness, is called the $k$-th osculating space. Here, the $k$-th osculating space is considered as a point of the appropriate Grassmannian $Gr\bigl(\binom{n+k}{k}-1,N\bigr)$ of $\binom{n+k}{k}-1$ dimensional projective subspaces of $\P^N$. In particular, the case $k=1$ gives the classical Gauss map. The functions \textsf{gausskFiber} and \textsf{gausskImage} allow for the computation of the general fiber and the image of the Gauss map of order $k$. To illustrate these methods we compute the fiber and image of the Gauss map of order $2$ for the singular toric variety associated to the polytope on the right using the code on the left.\newline 

\begin{minipage}{0.58\linewidth}
\begin{verbatim}
P=convexHull(transpose matrix{{0,0},{1,0},{3,1},{0,2}});
gausskFiber(latticePoints(P),2)
           2
= {1, x , x }
       1   1
gausskImage(latticePoints(P),2)
    4 5   5 5   6 5   7 5
= {x x , x x , x x , x x }
    0 1   0 1   0 1   0 1
\end{verbatim}
\end{minipage}
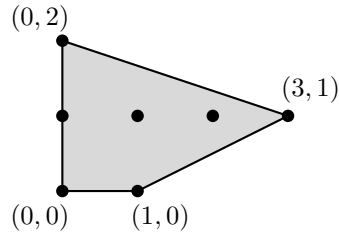
\begin{minipage}{0.2\linewidth}
\begin{tikzpicture}
\draw[draw=white] (4,3.5)--(5,3.5);
\draw[fill=black!15!white,thick, draw=black] (0,0)--(0,2)--(3,1)--(1,0)--(0,0);
\foreach \pt in {(0,0),(1,0),(0,2)}{
\node [fill=black,shape=circle, scale=0.5] at \pt {};
}
\foreach \y in {0,...,3}{
\node [fill=black,shape=circle, scale=0.5] at (\y,1) {};
}
\node at (-0.3,-0.35) {$(0,0)$};
\node at (1.3,-0.35) {$(1,0)$};
\node at (3.3,1.35) {$(3,1)$};
\node at (-0.3,2.3) {$(0,2)$};
\end{tikzpicture}

\end{minipage}\newline
The first command line constructs the polytope, while the second and third computes the fiber and image of the Gauss map of order $2$. 

For the remaining methods of the package \LP we refer to the documentation available through \textsf{Macaulay2}.

\subsection*{Acknowledgements}
We would like to thank Mats Boij for introducing us to Macaulay2 and proofreading. Moreover, we are greatful to Sandra Di Rocco and David Rydh for helpfull comments.
\bibliographystyle{plain}
\bibliography{LPart}

\begin{thebibliography}{10}

\bibitem{primer}
Thomas Bauer, Sandra Di~Rocco, Brian Harbourne, Micha{\l} Kapustka, Andreas
  Knutsen, Wioletta Syzdek, and Tomasz Szemberg.
\newblock A primer on {S}eshadri constants.
\newblock In {\em Interactions of classical and numerical algebraic geometry},
  volume 496 of {\em Contemp. Math.}, pages 33--70. AMS, Providence, RI, 2009.

\bibitem{Polyh}
Ren{\'e} Birkner.
\newblock Polyhedra: a package for computations with convex polyhedral objects.
\newblock {\em J. Softw. Algebra Geom.}, 1:11--15, 2009.

\bibitem{finiteSLP}
Tristram Bogart, Christian Haase, Milena Hering, Benjamin Lorenz, Benjamin
  Nill, Andreas Paffenholz, G{\"u}nter Rote, Francisco Santos, and Hal Schenck.
\newblock Finitely many smooth {$d$}-polytopes with {$n$} lattice points.
\newblock {\em Israel J. Math.}, 207(1):301--329, 2015.

\bibitem{C22}
Maria Castellani.
\newblock Sule superfici i cui spazi osculatori sono biosculatori.
\newblock {\em Rom. Acc. I. Rend.}, 5(31):347--350, 1922.

\bibitem{Cox}
David~A. Cox, John~B. Little, and Henry~K. Schenck.
\newblock {\em Toric varieties}, volume 124 of {\em Graduate Studies in
  Mathematics}.
\newblock AMS, Providence, RI, 2011.

\bibitem{SandraLinFib}
Sandra {Di Rocco}.
\newblock {Linear Toric Fibrations}.
\newblock In {\em Combinatorial Algebraic Geometry, Lecture Notes in
  Mathematics 2108}. Springer, 2014.

\bibitem{PolyAdj}
Sandra Di~Rocco, Christian Haase, Benjamin Nill, and Andreas Paffenholz.
\newblock Polyhedral adjunction theory.
\newblock {\em Algebra Number Theory}, 7(10):2417--2446, 2013.

\bibitem{HGM}
Sandra {Di Rocco}, Kelly {Jabbusch}, and Anders {Lundman}.
\newblock {A note on higher order Gauss maps}.
\newblock {\em ArXiv e-prints}, (1410.4811), 2014.

\bibitem{FukukawaIto}
Katsuhisa Furukawa and Atsushi Ito.
\newblock {Gauss maps of toric varieties}.
\newblock {\em ArXiv e-prints}, (1403.0793), 2014.

\bibitem{M2}
Daniel~R. Grayson and Michael~E. Stillman.
\newblock Macaulay2, a software system for research in algebraic geometry.
\newblock Available at \url{http://www.math.uiuc.edu/Macaulay2/}.

\bibitem{GH79}
Phillip Griffiths and Joseph Harris.
\newblock Algebraic geometry and local differential geometry.
\newblock {\em Ann. Sci. \'Ecole Norm. Sup. (4)}, 12(3):355--452, 1979.

\bibitem{degen}
Atsushi Ito.
\newblock {Seshadri constants via toric degenerations}.
\newblock {\em Journal f\"{u}r die reine und angewandte Mathematik},
  695:151--174, 2014.

\bibitem{Lazar}
Robert Lazarsfeld.
\newblock {\em Positivity in algebraic geometry. {I}}, volume~49 of {\em
  Ergebnisse der Mathematik und ihrer Grenzgebiete. 3. Folge. A Series of
  Modern Surveys in Mathematics}.
\newblock Springer-Verlag, Berlin, 2004.

\bibitem{lorenz}
Benjamin Lorenz.
\newblock {Classification of smooth lattice polytopes with few lattice points}.
\newblock Master thesis, Frei Universit\"{a}t Berlin, Fachbereich Mathematik.,
  2010.

\bibitem{myclass}
Anders Lundman.
\newblock A classification of smooth convex 3-polytopes with at most 16 lattice
  points.
\newblock {\em J. Algebraic Combin.}, 37(1):139--165, 2013.

\bibitem{localpos}
Anders Lundman.
\newblock Local positivity of line bundles on smooth toric varieties and
  {C}ayley polytopes.
\newblock {\em Journal of Symbolic Computation}, pages~--, 2015.

\bibitem{Perkinson}
David Perkinson.
\newblock {Inflections of Toric Varieties}.
\newblock {\em Michigan Mathematical Journal}, 48(1):483--515, 2000.

\bibitem{Pohl}
William~F. Pohl.
\newblock Differential geometry of higher order.
\newblock {\em Topology}, 1:169--211, 1962.

\end{thebibliography}

\end{document}